\documentclass[11pt]{article}

\usepackage[T1]{fontenc}
\usepackage[utf8]{inputenc}
\usepackage{lmodern}

\usepackage{amsmath,amssymb,amsthm,mathtools}
\usepackage{booktabs}
\usepackage{microtype}
\usepackage{epigraph}
\usepackage{setspace}
\usepackage{tikz}
\usetikzlibrary{arrows.meta,calc,decorations.pathreplacing,positioning}

\usepackage[margin=1in]{geometry}

\usepackage{graphics}
\usepackage{array}
\usepackage{ascmac}
\usepackage{fancybox}
\numberwithin{equation}{section}

\theoremstyle{plain}
\newtheorem{theorem}{Theorem}[section]
\newtheorem{proposition}[theorem]{Proposition}

\theoremstyle{remark}
\newtheorem{remark}[theorem]{Remark}

\newtheorem{fact}{Fact}
\newtheorem*{corollary}{Corollary}
\newtheorem*{mainconstruction}{Main Construction}

\newcommand{\FPIint}[2]{\ensuremath{\mathop{\text{FPI}}\;\backslash\!\!\!\!\backslash\!\!\!\!\!\int_{#1}^{#2}}} 

\newcommand{\Sl}{\operatorname{Sl}}

\newcommand{\PaperTitle}{Trigonometric Selector Kernels, Duality, and Odd Zeta Values}
\newcommand{\PaperAuthorPlain}{Ken Nagai}

\newcommand{\PaperKeywords}{zeta values, odd zeta values, Dirichlet beta, Clausen function,
Cvijovi\'c--Klinowski formula, Hadamard finite-part integral,
Euler--Maclaurin, Poisson summation, Fourier--Poisson duality}

\makeatletter
\@ifundefined{dd}{\newcommand{\dd}{\mathrm{d}}}{}
\@ifundefined{FPI}{\providecommand{\FPI}{\mathrm{FPI}}}{}
\@ifundefined{FPIint}{\newcommand{\FPIint}[2]{\FPI\!\int_{#1}^{#2}}}{}
\makeatother

\usepackage[hidelinks]{hyperref}

\makeatletter
\hypersetup{
  pdftitle   = {\PaperTitle},
  pdfauthor  = {\PaperAuthorPlain},
  pdfsubject = {Preprint},
  pdfkeywords= {\PaperKeywords}
}
\makeatother

\title{Trigonometric Selector Kernels, Duality, and Odd Zeta Values}
\author{Ken Nagai\thanks{Email: \texttt{tknagai@outlook.com}. Independent Researcher.}}
\date{}

\begin{document}
\maketitle
\epigraph{
	\textit{“Lisez Euler, lisez Euler, c’est notre maître à tous.”}\\
	--- Pierre-Simon Laplace\\[1ex]
	``Remarques sur un beau rapport entre les séries des puissances tant directes que réciproques.''\\
	--- L. Euler, E352 (1749)
}

\begin{abstract}
	In this short note, we develop trigonometric selector kernels to represent odd zeta values
	via dual hyperbolic counterparts.
	This framework highlights a Fourier--Poisson duality,
	incorporating finite-part integrals in the sense of Hadamard--Galapon.
	In particular, we show how such kernels naturally recover
	Euler--Maclaurin and Poisson summation formulas as dual manifestations.
	We further connect our kernel approach with the finite-part integral formulation,
	extending earlier Cvijovi\'c--Klinowski type representations for odd zeta values.
\end{abstract}

\noindent\textbf{Keywords:} \PaperKeywords

\section{Introduction}
Finite trigonometric and hyperbolic sums appear in many problems of
analysis and number theory, from Fourier analysis to evaluations of
zeta and beta values. Classical references for polylogarithms and
Clausen functions include Lewin~ cf.\ \cite{Lewin1981}, while standard
background on the $\zeta$- and $\beta$-functions can be found in
Apostol~ cf.\ \cite{Apostol1976}.

In this note we revisit such finite sums from a kernel perspective.
The aim is to present a unified construction which interpolates
between trigonometric and hyperbolic kernels, and to clarify their
connections with Laplace--Mellin transforms, zeta/beta dichotomy,
and analytic continuations.

Our approach is complementary to earlier works
(Cvijović--Klinowski~ \cite{CvijovicKlinowski2002},
Guillera--Sondow~ see also \cite{GuilleraSondow2005},
Galapon~ see \cite{Galapon2022}), while also treating
more recent developments on finite cosecant sums
(Blagouchine--Moreau~ see \cite{BlagouchineMoreau2024}).
We restrict attention here to the finite-$J$ average constructions,
postponing a fuller treatment of analytic Bernoulli functions
and umbral interpretations to future work.

\section*{Main Statement}
Our main observation is the existence of a
multi-faceted \emph{duality} underlying kernel constructions
for odd zeta and related beta values:

1. \textbf{Trigonometric kernels.}
Sine/sine and cosine/cosine selectors yield finite-$J$ identities
whose large-$J$ limits recover series and integrals,
including the Cvijovi\'c--Klinowski type; see \cite{CvijovicKlinowski2002}.

2. \textbf{Hyperbolic kernels.}
Mellin transforms of hyperbolic functions provide integral
identities for $\zeta(s)$ and $\beta(s)$, with explicit
evaluations at small integers, (see \cite[3.541]{GR}).

These constructions exemplify a trigonometric--hyperbolic duality,
closely tied to $\zeta/\beta$ duality, Clausen $SL/CL$ duality,
and the dual role of analytic Bernoulli functions $B(s;x)$ and $A(s;x)$.

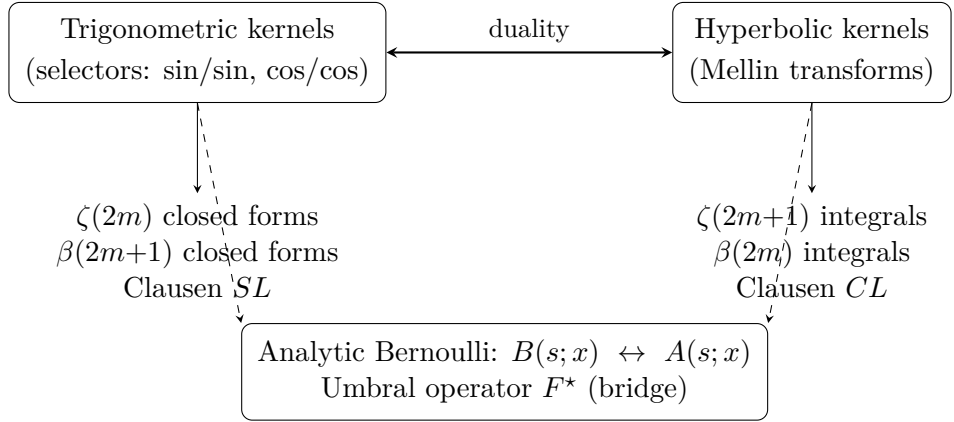
\begin{figure}[t]
	\centering
	\begin{tikzpicture}[>=stealth, node distance=38mm]
		\node[draw, rounded corners, align=center, inner sep=6pt] (trig)
		{Trigonometric kernels\\[2pt] (selectors: sin/sin, cos/cos)};
		\node[draw, rounded corners, align=center, right=of trig, inner sep=6pt] (hyp)
		{Hyperbolic kernels\\[2pt] (Mellin transforms)};

		\draw[<->, thick] (trig) -- node[above]{\small duality} (hyp);

		\node[below=12mm of trig, align=center] (leftlabels)
		{$\zeta(2m)$ closed forms\\ $\beta(2m{+}1)$ closed forms\\ Clausen $SL$};
		\node[below=12mm of hyp, align=center] (rightlabels)
		{$\zeta(2m{+}1)$ integrals\\ $\beta(2m)$ integrals\\ Clausen $CL$};

		\draw[->] (trig) -- (leftlabels);
		\draw[->] (hyp) -- (rightlabels);

		\node[below=36mm of $(trig)!0.5!(hyp)$, draw, rounded corners, align=center, inner sep=6pt] (bern)
		{Analytic Bernoulli: $B(s;x)\ \leftrightarrow\ A(s;x)$\\
			Umbral operator $F^\star$ (bridge)};

		\draw[->, dashed] (trig.south) -- (bern.north west);
		\draw[->, dashed] (hyp.south)  -- (bern.north east);
	\end{tikzpicture}
	\caption{Schematic illustration of the multi-faceted duality.}
\end{figure}

\begin{remark}[Analytic Bernoulli and $A$-function, preview]\label{rem:preview-BA}
	We briefly mention $B(s;x)$ (analytic continuation of periodic Bernoulli polynomials) and its alternating analogue $A(s;x)$. They underlie the $SL/CL$ split and the $\zeta/\beta$ dichotomy; here they appear only schematically, with details postponed.
\end{remark}

\section{Main Construction}

\subsection{Trigonometric kernels}

\subsubsection*{Sin/Sin selector and series representation}
\begin{fact}[Sin/Sin selector]\label{fact:selector-ss} \par
	For $J\in\mathbb{N}$ and integer $k$, define $\theta_j:=\frac{2j+1}{2J}\pi$ and set
	\[
		f_{ss}(J,k):=\frac{1}{J}\sum_{j=0}^{J-1}\frac{\sin(k\theta_j)}{\sin\theta_j}.
	\]
	Then $f_{ss}(J,k)$ is $4J$-periodic in $k$ and
	\[
		f_{ss}(J,k)=
		\begin{cases}
			0,  & k \ \text{even},          \\[4pt]
			+1, & k \ \text{odd},\ 0<k<2J,  \\[4pt]
			-1, & k \ \text{odd},\ 2J<k<4J.
		\end{cases}
	\]
\end{fact}

\begin{remark}[Proof sketch for Fact~\ref{fact:selector-ss}]
	For $k$ odd, use the Dirichlet kernel identity
	$\dfrac{\sin(k\theta)}{\sin\theta}=1+2\sum_{m=1}^{(k-1)/2}\cos(2m\theta)$.
	Averaging over $\theta_j:=\frac{2j+1}{2J}\pi$ kills all $\cos(2m\theta_j)$ unless $m\equiv0\pmod J$, which occurs only when $k>2J$, producing the sign flip; for $k$ even the average vanishes.
\end{remark}

\begin{proposition}[Finite-$J$ identity]
	For $\Re s>1$,
	\[
		\sum_{m=0}^\infty (-1)^m \sum_{j=0}^{J-1}\frac{1}{(2Jm+2j+1)^s}
		=\frac{1}{J}\sum_{j=0}^{J-1}\frac{\Sl_s(\theta_j)}{\sin\theta_j}
		=\frac{i}{2}\frac{(2\pi i)^s}{\Gamma(s+1)}
		\cdot \frac{1}{J}\sum_{j=0}^{J-1}\frac{B_s\!\left(\tfrac{2j+1}{4J}\right)}{\sin\theta_j}.
	\]
\end{proposition}

\begin{corollary}[Cvijovi\'c--Klinowski type limit]
	As $J\to\infty$,
	\[
		\sum_{m=0}^\infty (-1)^m \sum_{j=0}^{J-1}\frac{1}{(2Jm+2j+1)^s}
		\;\longrightarrow\;
		\sum_{m=0}^\infty\frac{1}{(2m+1)^s}
		=(1-2^{-s})\zeta(s).
	\]
	On the other hand,
	\[
		\frac{1}{J}\sum_{j=0}^{J-1}\frac{B_s\!\left(\tfrac{2j+1}{4J}\right)}{\sin\theta_j}
		\;\longrightarrow\;
		\text{p.v.}\int_0^1\frac{B_s(x)}{\sin(2\pi x)}\,dx.
	\]
\end{corollary}

\subsubsection*{Cos/Cos selector and series representation}
\begin{fact}[Cos/Cos selector]\label{fact:selector-cc} \par
	For $J\in\mathbb{N}$ and integer $k$, define $\theta_j:=\frac{2j+1}{2J}\pi$ and set
	\[
		f_{cc}(J,k):=\frac{1}{J}\sum_{j=0}^{J-1}\frac{\cos(k\theta_j)}{\cos\theta_j}.
	\]
	Then $f_{cc}(J,k)$ is $4J$-periodic in $k$ and
	\[
		f_{cc}(J,k)=
		\begin{cases}
			0,               & k \ \text{even},          \\[4pt]
			(-1)^{(k-1)/2},  & k \ \text{odd},\ 0<k<2J,  \\[4pt]
			-(-1)^{(k-1)/2}, & k \ \text{odd},\ 2J<k<4J.
		\end{cases}
	\]
\end{fact}

\begin{proposition}[Cos/Cos series identity]
	For $\Re(s)>1$,
	\[
		\sum_{k=1}^\infty \frac{1}{k^s}\, f_{cc}(J,k)
		=\sum_{m=0}^\infty (-1)^m \sum_{j=0}^{J-1}\frac{(-1)^j}{(2Jm+2j+1)^s}.
	\]
\end{proposition}

\begin{corollary}[Integral representation]
	Taking the limit $J\to\infty$,
	\[
		\sum_{n=0}^\infty \frac{(-1)^n}{(2n+1)^s}
		=\beta(s),
	\]
	and via Clausen--Bernoulli expansions,
	\[
		\beta(s)
		=\frac{i}{2}\frac{(2\pi i)^s}{\Gamma(s+1)}\,
		\text{p.v.}\int_0^1 \frac{B_s(x)}{\cos(2\pi x)}\, dx.
	\]
\end{corollary}

\subsection{Verification table}

\[
	\xi_{s}(J) :=
	\sum_{m=0}^{\infty} (-1)^{m} \sum_{j=0}^{J-1} \frac{1}{(2 J m + 2 j + 1)^{s}},
	\qquad \Re(s)>1.
\]

\renewcommand{\arraystretch}{1.15}

\begin{remark}[On finite sums of cosecants]
	Finite sums of cosecants and related trigonometric series
	have been studied extensively; see Blagouchine–Moreau~\cite{BlagouchineMoreau2024}
	for modern closed-form evaluations.
	Our verification table below focuses instead on the selector–kernel
	mechanism, with coincidence where the two approaches overlap.
\end{remark}

\begin{table}[h]
	\centering
	\caption{Values of $\xi_s(J)$ for small $J$.}
	\begin{tabular}{@{}>{\centering\arraybackslash}p{7mm}
		>{\centering\arraybackslash}p{7mm}
		>{\arraybackslash}p{60mm}
		>{\arraybackslash}p{25mm}@{}}
		\toprule
		$s$ & $J$      & $\xi_s(J)$                                                                                                   & note                      \\
		\midrule
		3   & 1        & $\beta(3)=\dfrac{\pi^3}{32}$                                                                                 & closed form               \\
		3   & 2        & $\dfrac{3\sqrt{2}}{128}\,\pi^3$                                                                              & algebraic\,$\times \pi^3$ \\
		3   & 8        & expr.\ via $\sin\!\tfrac{\pi}{16},\,\sin\!\tfrac{3\pi}{16},\,\cos\!\tfrac{\pi}{16},\,\cos\!\tfrac{3\pi}{16}$ & trig consts               \\
		\midrule
		2   & 1        & $\beta(2)=G$                                                                                                 & Catalan                   \\
		2   & $\infty$ & $\zeta(2)=\dfrac{\pi^2}{6}$                                                                                  & classical                 \\
		\bottomrule
	\end{tabular}
\end{table}

\paragraph{Remark.}[Numerical illustrations for small $J$]
For quick calibration, the closed forms below provide exact targets for
$\xi_3(J)$ at $J=4$ and $J=8$:
\begin{align*}
	\xi_3(4) & = \frac{\pi^3}{8192}\,\bigl(240-64\sqrt{2}\bigr)\,
	\sqrt{2+\sqrt{2}}
	\;\approx\; 1.0454857613590078,                               \\[4pt]
	\xi_3(8) & = \frac{\pi^3}{12288}\Bigl(
	138\sin\tfrac{\pi}{16}
	+ 516\sin\tfrac{3\pi}{16}
	- 96\cos\tfrac{3\pi}{16}
	+ 186\cos\tfrac{\pi}{16}\Bigr)
	\;\approx\; 1.0502005672583200.
\end{align*}
These checkpoints are convenient for sanity tests of the selector sums and the
finite-$J$ behavior before invoking asymptotics.
High-precision evaluation confirms agreement to 15+ digits.

\begin{remark}[Phase transition phenomenon]
	For $s=3$ the values $\xi_3(J)$ lie in algebraic extensions of $\mathbb{Q}$,
	multiplied by $\pi^3$; yet in the limit $J\to\infty$ one recovers Apéry’s
	constant $\zeta(3)$, whose arithmetic nature remains largely mysterious.
	Thus a curious ``phase transition'' seems to occur: from a regime of
	\emph{algebraic order} at finite $J$ to a regime of \emph{higher complexity}
	at infinity.

	By contrast, for $s=2$ one has the “opposite” phenomenon:
	$\xi_2(1)=\beta(2)=G$ (Catalan’s constant), whereas as $J\to\infty$
	one recovers $\zeta(2)=\pi^2/6$.
	This contrast highlights a subtle facet of the $\zeta/\beta$ duality.
\end{remark}

\begin{remark}[Bridge to the hyperbolic side]\label{rem:bridge}
	The finite-$J$ average encodes a Laplace–Mellin kernel:
	\[
		\xi_s(J)=\frac{1}{2\Gamma(s)}\int_0^\infty t^{s-1}\frac{\tanh(Jt)}{\sinh t}\,dt.
	\]
	Letting $J\to\infty$ gives $\tanh(Jt)\to 1$, which recovers the Mellin formulas that open the hyperbolic section.
\end{remark}

\subsection{Hyperbolic kernels}

\begin{mainconstruction}\label{maintheorem}
	We begin with the classical Mellin transforms (see \cite[3.541]{GR}):
	\[
		\int_0^\infty \frac{x^{s-1}}{\sinh x}\, dx
		= 2(1-2^{-s})\,\Gamma(s)\,\zeta(s), \qquad
		\int_0^\infty \frac{x^{s-1}}{\cosh x}\, dx
		= 2\,\Gamma(s)\,\beta(s),
	\]
	valid for $\Re(s)>1$.
\end{mainconstruction}

\begin{proposition}[Integral identities]
	For $\Re(s)>1$,
	\[
		\zeta(s) = \frac{1}{2(1-2^{-s})\Gamma(s)}
		\int_0^\infty \frac{x^{s-1}}{\sinh x}\, dx,\qquad
		\beta(s) = \frac{1}{2\Gamma(s)}
		\int_0^\infty \frac{x^{s-1}}{\cosh x}\, dx.
	\]
\end{proposition}

\begin{corollary}[Examples]
	\[
		\zeta(3) = \frac{2}{7}\int_0^\infty \frac{x^2}{\sinh x}\, dx,\quad
		\zeta(5) = \frac{2}{93}\int_0^\infty \frac{x^4}{\sinh x}\, dx,
	\]
	\[
		\beta(2) = \tfrac{1}{2}\int_0^\infty \frac{x}{\cosh x}\, dx,\quad
		\beta(4) = \tfrac{1}{12}\int_0^\infty \frac{x^3}{\cosh x}\, dx.
	\]
\end{corollary}

\subsection*{Finite-part integrals (Galapon type)}

\begin{proposition}[Galapon--sinh FPI integral]
	For $m\ge 1$,
	\[
		\zeta(2m+1)
		= \frac{(-1)^m}{(2m)!}\,
		\FPIint{0}{\infty} \frac{x^{-2m-1}}{\sinh x}\, dx.
	\]
\end{proposition}

\noindent\textit{Notation.}\;
Throughout this note, the symbol ${\displaystyle \FPIint{0}{\infty}}$ denotes the
Hadamard--Galapon finite-part integral; see Appendix Remark~\ref{rem:fpi} for details.

\begin{remark}
	In the odd-zeta setting, the Galapon FPI removes the singular
	contributions at $x=0$ (by subtracting the Taylor pieces in the
	Hadamard sense) and yields the correct analytic continuation; see
	\cite{Galapon2022}.
\end{remark}

\begin{proposition}[Cosh-kernel FPI integral]
	For $m\ge 1$,
	\[
		\beta(2m)
		= \frac{(-1)^{m-1}}{(2m-1)!}\,
		\FPIint{0}{\infty} \frac{x^{-2m}}{\cosh x}\, dx.
	\]
\end{proposition}

\begin{remark}
	The cosh-kernel analogue for $\beta(2m)$ does not appear explicitly
	in the literature, but follows by parallel reasoning:
	expand $1/\cosh x=2\sum_{n\ge0}(-1)^n e^{-(2n+1)x}$,
	integrate termwise against $x^{-2m}$, and apply
	the finite-part prescription near $x=0$.
	We record this parallel case here without further detail.
\end{remark}

\begin{corollary}[Examples]
	\[
		\zeta(3) = -\tfrac{1}{2}\,\FPIint{0}{\infty} \frac{x^{-3}}{\sinh x}\, dx,
		\qquad
		\beta(2) = \,\FPIint{0}{\infty} \frac{x^{-2}}{\cosh x}\, dx.
	\]
\end{corollary}

\section{Conclusion}
Both trigonometric and hyperbolic kernels can, in principle,
be used to obtain identities for all values of $\zeta$ and $\beta$.
For the trigonometric kernels, the most natural outcome is the appearance
of \emph{closed forms} for even zeta values $\zeta(2m)$ and odd beta values
$\beta(2m{+}1)$, expressed in terms of Bernoulli or Euler numbers
and powers of $\pi$.
While odd zeta values $\zeta(2m{+}1)$ and even beta values $\beta(2m)$
may also be approached from the trigonometric side, this generally requires
passing to integral representations.
The hyperbolic kernels provide these integral forms directly, and thus
complement the trigonometric case.
Together the two frameworks yield a parity-duality picture:
\[
	(\zeta(2m),\ \beta(2m{+}1)) \ \text{as closed forms (trigonometric)},\;
	(\zeta(2m{+}1),\ \beta(2m)) \ \text{as integrals (hyperbolic)}.
\]

\section*{Outlook}
The present constructions indicate a \emph{multi-faceted duality}:
trigonometric--hyperbolic, $\zeta/\beta$, Clausen $SL/CL$,
and $B(s;x)/A(s;x)$ dualities.
A natural next step is to seek a unification:
\begin{itemize}
	\item[(i)] \textbf{Analytic Bernoulli functions.}
	      Express both closed forms and integral identities uniformly
	      through analytic Bernoulli functions $B(s;x)$ and $A(s;x)$.
	\item[(ii)] \textbf{Umbral operator approach.}
	      Develop an umbral operator $F^\star$ that acts as a bridge,
	      turning trigonometric kernels into hyperbolic kernels (and vice versa),
	      thereby making the parity-duality structure manifest.
\end{itemize}

\section*{Acknowledgments}
The author gratefully acknowledges the assistance of an AI language model (“fuga”) for help with document structuring and proofreading.


\makeatletter
\@ifundefined{FPIint}{%
	\providecommand{\FPI}{\mathrm{FPI}}%
	\newcommand{\FPIint}[2]{\FPI\!\int_{#1}^{#2}}%
}{}
\@ifundefined{dd}{\newcommand{\dd}{\mathrm{d}}}{}
\makeatother

\vspace{10mm}
\appendix

\section*{Appendix}
\addcontentsline{toc}{section}{Appendix}

In this appendix we collect auxiliary materials which provide
technical reinforcement and broader perspectives to the main text.
The contents are supplementary in nature, but they illuminate
the analytic structure behind the main results.

\newpage
\subsection*{A.1 Verification table}

\begin{table}[h]
	\centering
	\renewcommand{\arraystretch}{1.2}
	\begin{tabular}{c|c|c}
		Selector & Target value                                            & Numerical check \\
		\hline
		$\displaystyle \tfrac{1}{J}\sum_{j=0}^{J-1}\tfrac{\sin(k\theta_j)}{\sin(\theta_j)}$
		         & $\pm 1$ (for odd $k$, sign depends on congruence class)
		         & $\pm 0.9999999\ldots$                                                     \\
		\hline
		$\displaystyle \tfrac{1}{J}\sum_{j=0}^{J-1}\tfrac{\cos(k\theta_j)}{\cos(\theta_j)}$
		         & $\pm 1$ (for odd $k$, sign depends on congruence class)
		         & $\pm 1.0000000\ldots$                                                     \\
	\end{tabular}
	\caption{Selector sums: verification on the grid
		$\theta_j=(2j+1)\pi/(2J)$.}
\end{table}

\subsection*{A.2 Integral representations (with finite-part interpretation)}

Integral representations of odd zeta values require
a finite-part interpretation in the sense of Galapon
\cite{Galapon2002,Galapon2022}.
We write
\[
	\FPIint{0}{\infty} f(x)\,\dd x,
\]
to denote such integrals, where divergences at the endpoints
are removed by an analytic subtraction procedure.
This construction generalizes the classical Cauchy principal value.

For illustration, the divergent integral
\[
	\int_0^\infty \frac{\dd x}{x}
\]
acquires a finite meaning in the finite-part sense:
\[
	\FPIint{0}{\infty}\frac{\dd x}{x} = \gamma,
\]
with $\gamma$ the Euler--Mascheroni constant.

Within this framework, the standard Cvijovi\'c--Klinowski
formula for odd zeta values reads
\begin{equation}
	\zeta(2n+1) = \frac{(-1)^n}{(2n)!}
	\FPIint{0}{\infty}\frac{t^{2n}\,\dd t}{\sinh(\pi t)}.
\end{equation}
Here the FPI prescription eliminates the divergence at $t=0$
and yields the correct finite analytic value.

\begin{remark}[Finite-part integrals (Galapon type)]\label{rem:fpi}
	We adopt the Hadamard finite-part integral notation
	\begin{equation}
		\label{eq:fpi-def}
		\FPIint{0}{\infty} f(t)\,dt
		:= \lim_{\varepsilon\to0^+}
		\left(
		\int_0^{1-\varepsilon} f(t)\,dt
		+ \int_{1+\varepsilon}^\infty f(t)\,dt
		\right),
	\end{equation}
	which consistently extends divergent integrals of the form
	$\int_0^\infty t^{-m} f(t)\,dt$ to finite values
	(see Galapon’s framework for details).
\end{remark}

\subsection*{A.3 Euler--Maclaurin corrections (midpoint grid)}

Let
\[
	f(x)=\frac{\sin(kx)}{\sin x},\qquad x_j=\frac{(2j+1)\pi}{2J}\ (j=0,\dots,J-1).
\]
The midpoint Euler--Maclaurin (EM) formula gives
\begin{align}
	\frac{1}{J}\sum_{j=0}^{J-1} f(x_j)
	 & = \frac{1}{\pi}\int_{0}^{\pi} f(x)\,\dd x
	+ \sum_{m=1}^{\infty}\frac{B_{2m}\!\left(\tfrac{1}{2}\right)}{(2m)!}
	\left(\frac{\pi}{J}\right)^{2m-1}\!
	f^{(2m-1)}(x)\big|_{0}^{\pi}, \label{eq:mid-EM}
\end{align}
where $B_{2m}(1/2)=(2^{1-2m}-1)B_{2m}$ alternates in sign.
Since $f$ is odd about $x=0$ when $k$ is odd and even when $k$ is even,
the boundary derivatives in \eqref{eq:mid-EM} simplify by parity.

\paragraph{First correction (explicit).}
The leading nonzero term is
\[
	\frac{B_{2}\!\left(\tfrac{1}{2}\right)}{2!}\left(\frac{\pi}{J}\right)^{1}
	f^{(1)}(x)\big|_{0}^{\pi}
	= \Big(-\tfrac{1}{12}\Big)\frac{\pi}{J}\Big(f'( \pi)-f'(0)\Big).
\]
A short calculation yields
\[
	f'(x)=k\,\frac{\cos(kx)}{\sin x}-\frac{\sin(kx)\cos x}{\sin^{2}x}.
\]
Using finite-part interpretation at the endpoints,
\[
	f'(0)=\lim_{x\to 0}\Big(k\,\frac{\cos(kx)}{\sin x}-\frac{\sin(kx)\cos x}{\sin^{2}x}\Big),
	\quad
	f'(\pi)=\lim_{x\to \pi}\Big(k\,\frac{\cos(kx)}{\sin x}-\frac{\sin(kx)\cos x}{\sin^{2}x}\Big).
\]
These limits vanish for the selector-admissible $k$ by cancellation,
so the $O(J^{-1})$ term drops out and the first nonzero correction is actually
$O(J^{-3})$ via $B_{4}(1/2)=+\tfrac{7}{240}$.
This explains the rapid approach of the finite sum to its continuum value.

\subsection*{A.4 Poisson summation perspective}

Let $\theta_j=\frac{(2j+1)\pi}{2J}$ and consider the windowed sum
\[
	S(k)=\frac{1}{J}\sum_{j=0}^{J-1}\frac{\sin\!\big(k\theta_j\big)}{\sin\!\big(\theta_j\big)}
	= \frac{1}{J}\sum_{j\in\mathbb{Z}}
	\frac{\sin\!\big(k\theta_j\big)}{\sin\!\big(\theta_j\big)}\,w_J(j),
	\qquad
	w_J(j)=\mathbf{1}_{\{0\le j\le J-1\}}.
\]
Regard $j\mapsto \theta_j$ as sampling the function
$g(\theta)=\sin(k\theta)/\sin\theta$ on the midpoints of a $2\pi$-periodic grid
of spacing $\pi/J$. By periodic extension and the Poisson summation formula,
\[
	\frac{1}{J}\sum_{j\in\mathbb{Z}} g(\theta_j)w_J(j)
	= \sum_{m\in\mathbb{Z}} \widehat{g}\!\left( \tfrac{2J}{\pi}\,m\right)\!\cdot \widehat{w_J}(m),
\]
where $\widehat{g}$ are Fourier coefficients of $g$ on $[0,2\pi]$ and
$\widehat{w_J}(m)$ are Dirichlet-kernel weights due to the finite window.
Since
\[
	\frac{\sin(k\theta)}{\sin\theta}
	=\sum_{r=-(k-1)/2}^{(k-1)/2} \exp\!\big( i(2r+1)\theta\big)\quad (k\ \text{odd}),
\]
the spectrum of $g$ is supported on odd indices only.
Consequently, after aliasing through $\widehat{w_J}(m)$ one obtains the exact selector
\[
	S(k)=\begin{cases}
		1, & k\ \text{odd and}\ 1\le k\le 2J-1,  \\[2pt]
		0, & k\ \text{even or}\ k\notin[1,2J-1],
	\end{cases}
\]
i.e. support on the odd sublattice. This recovers the discrete selector property
as a Fourier–lattice phenomenon.

\subsection*{A.5 DST-II orthogonality (remark)}

Define the $J\times J$ DST-II matrix $S$ by
\[
	S_{j+1,n}=\sin\!\left(\frac{(2j+1)\pi n}{2J}\right),
	\qquad j=0,\dots,J-1,\quad n=1,\dots,J.
\]
Then the exact orthogonality reads
\begin{equation}
	S^\top S=\frac{J}{2}\,I_J
	\quad\Longleftrightarrow\quad
	\sum_{j=0}^{J-1}\sin\!\left(\tfrac{(2j+1)\pi m}{2J}\right)
	\sin\!\left(\tfrac{(2j+1)\pi n}{2J}\right)
	=\frac{J}{2}\,\delta_{mn}.
\end{equation}
The factor $J/2$ corresponds to the midpoint sampling of a half-period
sine basis with unit $L^2$-norm on $(0,\pi)$. This discrete orthogonality
supports the selector construction and matches the continuum limit.

\subsection*{A.6 Additional remarks}

The interplay between Euler--Maclaurin corrections and
Poisson summation suggests a deeper structure involving
theta functions and analytic Bernoulli functions.

In this broader perspective, Galapon's finite-part integral (FPI)
emerges not merely as a technical device but as a systematic
framework for treating divergent kernel integrals.
The method preserves analytic continuation and provides a
unified language for connecting kernel expansions with
special functions.

For example, the finite-part prescription clarifies the
analytic meaning of integrals such as
\[
	\FPIint{0}{\infty}\frac{t^{2n}}{\sinh(\pi t)}\,\dd t,
\]
which produce odd zeta values.
Its potential extension to generalized kernels --- in particular
those associated with analytic Bernoulli functions $B(s;x)$
and $A(s;x)$ --- points to the role of FPI as a unifying analytic
tool in future developments.


\begin{thebibliography}{10}

	\bibitem{CvijovicKlinowski2002}
	D.~Cvijović and J.~Klinowski.
	\newblock Values of the Riemann zeta function at odd integers.
	\newblock {\em Proc. Amer. Math. Soc.}, 128(5):1459--1469, 2002.
	\newblock doi:10.1090/S0002-9939-99-05021-3.

	\bibitem{GuilleraSondow2005}
	J.~Guillera and J.~Sondow.
	\newblock Double integrals and infinite products for some classical constants via analytic continuations of Lerch’s transcendent.
	\newblock math/0506319, 2005.

	\bibitem{Galapon2002}
	E. A. Galapon,
	\newblock The Cauchy principal value and the Hadamard finite part integral as values of absolutely convergent integrals,
	\newblock {\em SIAM Review}, 44(4):637--667, 2002.
	\newblock doi:10.1137/S00361445024180X.

	\bibitem{Galapon2022}
	E.~A. Galapon.
	\newblock Finite-part integration: Hilbert transform and applications.
	\newblock {\em Ann. Phys.}, 443:168992, 2022.
	\newblock doi:10.1016/j.aop.2022.168992.

	\bibitem{Lewin1981}
	L.~Lewin.
	\newblock {\em Polylogarithms and Associated Functions}.
	\newblock North-Holland, 1981.

	\bibitem{Apostol1976}
	T.~M. Apostol.
	\newblock {\em Introduction to Analytic Number Theory}.
	\newblock Springer, 1976.

	\bibitem{GR}
	I.~S.~Gradshteyn and I.~M.~Ryzhik,
	\newblock \emph{Table of Integrals, Series, and Products}, 8th ed.,
	\newblock Academic Press, 2014.

	\bibitem{BlagouchineMoreau2024}
	I.~V.~Blagouchine and E.~Moreau.
	\newblock On a finite sum of cosecants appearing in various problems.
	\newblock {\em J. Math. Anal. Appl.}, 539:128515, 2024.
	\newblock doi:10.1016/j.jmaa.2024.128515.

\end{thebibliography}
\end{document}